\documentclass[11pt, final]{article}
\usepackage{a4}
\usepackage{amsmath}%
\usepackage{amstext}%
\usepackage{amssymb}%
\usepackage{showkeys}%
\usepackage{epsfig}%
\usepackage{cite}
\usepackage{tikz}
\usepackage{subfig}
\usepackage{caption}
\usepackage{multirow}
\usepackage{booktabs}
\usepackage{floatrow}

\usepackage{listings}
\usepackage[margin=0.5in]{geometry}
\newtheorem{theorem}{Theorem}

\newtheorem{algorithm}[theorem]{Algorithm}

\newtheorem{pb}[theorem]{Problem}
\newenvironment{problem}{\pb\rm}{\endpb}

\newtheorem{rem}[theorem]{Remark}
\newenvironment{remark}{\rem\rm}{\endrem}

\newcounter{unnumber}

\newtheorem{prf}[unnumber]{Proof.}
\newenvironment{proof}{\prf\rm}{\hfill{$\blacksquare$}\endprf}
\newcommand{\R}{\mathbb{R}}%
\newcommand{\N}{\mathbb{N}}%
\newcommand{\ol}{\overline}%

\DeclareMathOperator*\inte{int}%
\DeclareMathOperator*\sqri{sqri}%
\DeclareMathOperator*\ri{ri}%
\DeclareMathOperator*\dom{dom}%
\DeclareMathOperator*\B{\overline{\R}}%
\DeclareMathOperator*\gr{Gr}%
\DeclareMathOperator*\ran{ran}%
\DeclareMathOperator*\id{Id}%
\DeclareMathOperator*\prox{prox}%
\DeclareMathOperator*\argmin{argmin}

\DeclareMathOperator*\zer{zer}

\title{An inertial alternating direction method of multipliers}

\author{Radu Ioan Bo\c{t} \thanks{University of Vienna, Faculty of Mathematics, Oskar-Morgenstern-Platz 1, A-1090 Vienna, Austria,
email: radu.bot@univie.ac.at. Research partially supported by DFG (German Research Foundation), project BO 2516/4-1.} \and
Ern\"{o} Robert Csetnek \thanks {University of Vienna, Faculty of Mathematics, Oskar-Morgenstern-Platz 1, A-1090 Vienna, Austria,
email: ernoe.robert.csetnek@univie.ac.at. Research supported by DFG (German Research Foundation), project BO 2516/4-1.}}

\begin{document}
\maketitle

\noindent \textbf{Abstract.} In the context of convex optimization problems in Hilbert spaces, we induce inertial effects into the classical ADMM numerical scheme and obtain in this way so-called inertial ADMM algorithms, the convergence properties of which we investigate into detail. To this aim we make use of the inertial version of the Douglas-Rachford splitting method for monotone inclusion problems recently introduced in \cite{b-c-h-inertial}, in the context of concomitantly solving a convex minimization problem and its Fenchel dual. The convergence of both sequences of the generated iterates and of the objective function values is addressed. We also show how the obtained results can be extended to the treating of
convex minimization problems having as objective a finite sum of convex functions.
\vspace{1ex}

 \vspace{1ex}


\section{Introduction}\label{sec-intr}

One of the most popular algorithms in the literature for solving the convex optimization problem
\begin{equation}\label{prim-fin}  
\inf_{x\in\R^n}\{f(x)+g(Ax)\},
\end{equation}
where $f:\R^n\rightarrow\B$ and $g:\R^m\rightarrow\B$ are proper, convex and lower semicontinuous functions and $A$ is a $m \times n$ matrix with real entries, is the {\it alternating direction 
method of multipliers} (ADMM). We briefly describe this procedure. By introducing an auxiliary variable one can rewrite \eqref{prim-fin} as  
\begin{equation}\label{prim-fin-x-z}  \inf_{\substack{(x,z)\in\R^n\times\R^m\\Ax-z=0}}\{f(x)+g(z)\}.
\end{equation}
For $\gamma \geq 0$ we consider the {\it augmented Lagrangian} ${\cal L}_{\gamma}:\R^n\times\R^m\times\R^m\rightarrow\B$ defined by 
$${\cal L}_{\gamma}(x,z,y)=f(x)+g(z)+y^T(Ax-z)+\frac{\gamma}{2}\|Ax-z\|^2  \ \forall(x,y,z)\in\R^n\times\R^m\times\R^m,$$
where the Euclidean norm on $\R^m$ is taken. 

The ADMM algorithm reads for given $y^0, z^0 \in \R^m$ and every $k \geq 0$
\begin{equation}x^{k+1}=\argmin_{x\in\R^n} {\cal L}_{\gamma}(x,z^k,y^k) 
\end{equation}
\begin{equation}z^{k+1}=\argmin_{z\in\R^m} {\cal L}_{\gamma}(x^{k+1},z,y^k) 
\end{equation}
\begin{equation}y^{k+1}=y^k+\gamma(Lx^{k+1}-z^{k+1}).
\end{equation}

The convergence of the ADMM algorithm is guaranteed by assuming that the matrix $A$ has full column rank and the unaugmented Lagrangian $L_0$ has a {\it saddle point} 
$(\ol x,\ol z,\ol y)\in \R^n\times\R^m\times\R^m$, that is 
$${\cal L}_0(\ol x,\ol z,y)\leq {\cal L}_0(\ol x,\ol z,\ol y)\leq {\cal L}_0(x,z,\ol y) \ \forall (x,z,y)\in\R^n\times\R^m\times\R^m.$$
Let us mention that if $(\ol x,\ol z,\ol y)$ is a saddle point of ${\cal L}_0$, then $\ol x$ is an optimal solution to \eqref{prim-fin}, $\ol z=A\ol x$ and 
$\ol y$ is an optimal solution to the Fenchel dual problem to \eqref{prim-fin}
\begin{equation}\label{dual-fin}  \sup_{v\in\R^m}\{-f^*(-A^Tv)-g^*(v)\},
\end{equation}
where $A^T$ denotes the transpose of the matrix $A$ and $f^*$ and $g^*$ the conjugate functions of $f$ and $g$, respectively. 

One of the limitations of this algorithm is the presence of the term $Ax$ in the update rule of $x^{k+1}$, which means that the scheme is not a really full splitting algorithm, like the primal-dual algorithms recently considered in \cite{br-combettes, b-c-h1, b-c-h2}. Nevertheless, the algorithm has been successfully implemented in the context of different real-life problems, like location problems, the lasso problem in image processing, problems arising in satistics, support vector machines classification, etc. We refer 
the reader to the seminal work \cite{bpcpe} for the history of the ADMM algorithm and various concrete applications of it (see also \cite{ecb, eck, gabay, esser}). 

In this paper we propose new ADMM type numerical schemes, which have their roots in the class of so-called inertial proximal point algorithms.  The latter iterative schemes are designed for solving monotone inclusion problems and, as they arise from the time discretization of some differential inclusions of second order type (see \cite{alvarez2000, alvarez-attouch2001}), have the property that the next iterate is defined by using the previous two iterates. In this way an inertial effect is induced into the numerical scheme, the increasing interest in this class of algorithms being emphasized by a considerable number of papers written in the last fifteen years on this topic, see  \cite{alvarez2000, alvarez-attouch2001, alvarez2004, att-peyp-red, b-c-inertial, b-c-h-inertial, mainge2008, mainge-moudafi2008, moudafi-oliny2003, cabot-frankel2011}. 

We derive the inertial version of the ADMM from the perspective of the monotone operator theory, using as starting point the fact pointed out in \cite{gabay} that  the classical ADMM can be approached from the Douglas-Rachford splitting scheme for monotone inclusion problems (see also \cite{ecb}). In \cite{b-c-h-inertial} we recently introduced and studied the convergence properties of an {\it inertial Douglas-Rachford splitting algorithm}. By combining this iterative scheme with the techniques from \cite{gabay, ecb}, we are able to obtain an inertial ADMM scheme for simultaneously solving convex minimization problems 
and their Fenchel-type duals. For the sake of generality, the analysis is carried out in infinite dimensional Hilbert spaces, in opposition to the usual literature on  ADMM algorithms where the finite dimensional setting is preferred. Moreover, we prove the convergence of both sequences of the generated iterates and of the objective function values and show that the classical ADMM scheme can be recovered as particular instance of our inertial ADMM algorithm. We also point out how other ADMM-type algorithms from the literature turn out to be particular schemes of the new ones presented here. 

The paper is organized as follows. In the next section we make the reader familiar with the notions and results which will be  used throughout the manuscript. In Section \ref{IADMM} we introduce the inertial ADMM algorithm for simultaneously solving in Hilbert spaces the convex optimization problems which assumes the minimization of the sum of a proper, convex and lower semicontinuous function with the composition of another proper, convex and lower semicontinuous function with a linear continuous operator and its Fenchel dual problem and study its convergence  properties. Finally, in the last section we treat the 
convex minimization problem having as objective the finite sum of proper, convex and lower semicontinuous functions and its Fenchel-type dual and provide for this primal-dual pair inertial ADMM algorithms and corresponding convergence statements.

\section{Preliminaries}\label{prel}

For the readers convenience let us recall some standard notions and results in monotone operator theory and convex analysis which will be used further in the paper, see also 
\cite{bo-van, b-hab, bauschke-book, EkTem, simons, Zal-carte}. Let $\N= \{0,1,2,...\}$ be the set of nonnegative integers.
Let ${\cal H}$ be a real Hilbert space with \textit{inner product} $\langle\cdot,\cdot\rangle$ and associated \textit{norm} $\|\cdot\|=\sqrt{\langle \cdot,\cdot\rangle}$.
The symbols $\rightharpoonup$ and $\rightarrow$ denote weak and strong convergence, respectively.
When ${\cal G}$ is another Hilbert space and $L:{\cal H} \rightarrow {\cal G}$ a linear continuous operator,
then $L^* : {\cal G} \rightarrow {\cal H}$, defined by $\langle L^*y,x\rangle = \langle y,Lx \rangle$ for all
$(x,y) \in {\cal H} \times {\cal G}$, denotes the \textit{adjoint operator} of $L$.

For an arbitrary set-valued operator $A:{\cal H}\rightrightarrows {\cal H}$ we denote by 
$\gr A=\{(x,u)\in {\cal H}\times {\cal H}:u\in Ax\}$ its \emph{graph} and by $A^{-1}:{\cal H}\rightrightarrows {\cal H}$ its
\emph{inverse operator}, defined by $(u,x)\in\gr A^{-1}$ if and only if $(x,u)\in\gr A$.
We use also the notation $\zer A=\{x\in{\cal{H}}:0\in Ax\}$ for the \emph{set of zeros} of $A$. We say that $A$ is \emph{monotone}
if $\langle x-y,u-v\rangle\geq 0$ for all $(x,u),(y,v)\in\gr A$. A monotone operator $A$ is said to be \emph{maximally monotone}, if there exists no proper monotone extension of the graph of $A$ on ${\cal H}\times {\cal H}$.
The \emph{resolvent} of $A$, $J_A:{\cal H} \rightrightarrows {\cal H}$, is defined by $J_A=(\id_{{\cal H}}+A)^{-1}$, where $\id_{{\cal H}} :{\cal H} \rightarrow {\cal H}, \id_{\cal H}(x) = x$ for all $x \in {\cal H}$, is the \textit{identity operator} on ${\cal H}$. Moreover, if $A$ is maximally monotone, then $J_A:{\cal H} \rightarrow {\cal H}$ is single-valued and maximally monotone
(see \cite[Proposition 23.7 and Corollary 23.10]{bauschke-book}). For an arbitrary $\gamma>0$ we have (see \cite[Proposition 23.2]{bauschke-book})
\begin{equation}p\in J_{\gamma A}x \ \mbox{if and only if} \ (p,\gamma^{-1}(x-p))\in\gr A.\end{equation}

The operator $A$ is said to be \textit{uniformly monotone} if there exists an increasing function
$\phi_A : [0,+\infty) \rightarrow [0,+\infty]$ that vanishes only at $0$, and
$\langle x-y,u-v \rangle \geq \phi_A \left( \| x-y \|\right)$ for every $(x,u)\in\gr A$ and $(y,v) \in \gr A$. A well-known 
class of operators fulfilling this property is the one of the strongly monotone operators.
Let $\gamma>0$ be arbitrary. We say that  $A$ is \textit{$\gamma$-strongly monotone}, 
if $\langle x-y,u-v\rangle\geq \gamma\|x-y\|^2$ for all $(x,u),(y,v)\in\gr A$.

Let us recall now some elements of convex analysis. For a function $f:{\cal H}\rightarrow\overline{\R}$, where $\overline{\R}:=\R\cup\{\pm\infty\}$ is the extended real line, 
we denote by $\dom f=\{x\in {\cal H}:f(x)<+\infty\}$ its \textit{effective domain} and say that $f$ is \textit{proper} if $\dom f\neq\varnothing$ and $f(x)\neq-\infty$ for all $x\in {\cal H}$. 
We denote by $\Gamma({\cal H})$ the family of proper, convex and lower semi-continuous extended real-valued functions defined on ${\cal H}$. 
Let $f^*:{\cal H} \rightarrow \overline \R$, $f^*(u)=\sup_{x\in {\cal H}}\{\langle u,x\rangle-f(x)\}$ for all $u\in {\cal H}$, be the \textit{conjugate function} of $f$. 
The \textit{subdifferential} of $f$ at $x\in {\cal H}$, with $f(x)\in\R$, is the set $\partial f(x):=\{v\in {\cal H}:f(y)\geq f(x)+\langle v,y-x\rangle \ \forall y\in {\cal H}\}$. 
We take by convention $\partial f(x):=\varnothing$, if $f(x)\in\{\pm\infty\}$.  Notice that if $f\in\Gamma({\cal H})$, then $\partial f$ is a maximally monotone operator (see \cite{rock}) and it holds $(\partial f)^{-1} =
\partial f^*$. Let $S\subseteq {\cal H}$ be a nonempty set. The \textit{indicator function} of $S$, $\delta_S:{\cal H}\rightarrow \overline{\R}$, 
is the function which takes the value $0$ on $S$ and $+\infty$ otherwise. The subdifferential of the indicator function is the \textit{normal cone} 
of $S$, that is $N_S(x)=\{u\in {\cal H}:\langle u,y-x\rangle\leq 0 \ \forall y\in S\}$, if $x\in S$ and $N_S(x)=\emptyset$ for $x\notin S$. Notice that, if $S$ is a linear subspace, then $N_S(x)=S^{\perp}=\{u\in{\cal H}:\langle y,u\rangle=0 \ \forall y\in S\}$ for all $x\in S$. 

When $f\in\Gamma({\cal H})$ and $\gamma > 0$, for every $x \in {\cal H}$ we denote by $\prox_{\gamma f}(x)$ the \textit{proximal point} of parameter $\gamma$ of $f$ at $x$, 
which is the unique optimal solution of the optimization problem
\begin{equation}\label{prox-def}\inf_{y\in {\cal H}}\left \{f(y)+\frac{1}{2\gamma}\|y-x\|^2\right\}.
\end{equation}

Notice that the resolvent of the maximally monotone operator $\partial f$ is nothing else than the proximal point operator of $f$, namely,
\begin{equation}J_{\gamma\partial f}=(\id\nolimits_{\cal H}+\gamma\partial f)^{-1}=\prox\nolimits_{\gamma f}.\end{equation}
Moreover, if $f=\delta_S$, where $S\subseteq{\cal H}$ is a nonempty, closed convex set, then the proximal point operator of $f$ is the \textit{orthogonal projection} on $S$. 

Let us also recall that a proper function $f:{\cal H} \rightarrow \overline \R$ is said to be \textit{uniformly convex}, if there exists
an increasing function $\phi :[0,+\infty)\rightarrow[0,+\infty]$ which vanishes only at $0$ and such that
$$f(t x+(1-t)y)+t(1-t)\phi(\|x-y\|)\leq tf(x)+(1-t)f(y) \ \forall x,y\in\dom f\mbox{ and } \forall t\in(0,1).$$
In case this inequality holds for $\phi=(\beta/2)(\cdot)^2$, where $\beta >0$, then $f$ is said to be
\textit{$\beta$-strongly convex}. Let us mention that this property implies $\beta$-strong monotonicity of $\partial f$ (see \cite[Example 22.3]{bauschke-book})
(more general, if $f$ is uniformly convex, then $\partial f$ is uniformly monotone, see \cite[Example 22.3]{bauschke-book}).

We close this section by presenting the inertial Douglas-Rachford splitting algorithm for determining the zeros of the sum of two maximally monotone operators recently obtained in  \cite{b-c-h-inertial}, which will be crucial for the proof of the main results in the next section.

\begin{theorem}\label{inertial-DR} (Inertial Douglas--Rachford splitting algorithm, see \cite{b-c-h-inertial}) 
Let $A,B:{\cal H}\rightrightarrows{\cal H}$ be 
maximally monotone operators such that $\zer(A+B)\neq\varnothing$. Consider the following iterative scheme: 
$$(\forall k\geq 1)\hspace{0.2cm}\left\{
\begin{array}{ll}
y^k=J_{\gamma B}[w^k+\alpha_k(w^k-w^{k-1})]\\
v^k=J_{\gamma A}[2y^k-w^k-\alpha_k(w^k-w^{k-1})]\\
w^{k+1}=w^k+\alpha_k(w^k-w^{k-1})+\lambda_k(v^k-y^k)
\end{array}\right.$$
where $\gamma>0$, $w^0,w^1$ are arbitrarily chosen in $\cal{H}$, $(\alpha_k)_{k\geq 1}$ is nondecreasing with  $\alpha_1=0$ and 
$0\leq \alpha_k\leq\alpha < 1$ for every $k\geq 1$ and $\lambda, \sigma, \delta >0$ are such that 
\begin{equation}\label{delta-lambda}
\delta>\frac{\alpha^2(1+\alpha)+\alpha\sigma}{1-\alpha^2} \ \mbox{and} \ 0<\lambda\leq\lambda_k\leq 2\cdot\frac{\delta-\alpha\Big[\alpha(1+\alpha)+\alpha\delta+\sigma\Big]}{\delta\Big[1+\alpha(1+\alpha)+\alpha\delta+\sigma\Big]} \
\forall k \geq 1.
\end{equation}
Then there exists $x\in{\cal H}$ such that 
the following statements are true:
\begin{itemize} \item[(i)] $J_{\gamma B}x\in\zer(A+B)$;
\item[(ii)] $\sum_{k\in\N}\|w^{k+1}-w^k\|^2<+\infty$; 
\item[(iii)] $(w^k)_{k\in\N}$ converges weakly to $x$;
\item[(iv)] $y^k-v^k\rightarrow 0$ as $k\rightarrow+\infty$;
\item[(v)] $(y^k)_{k\geq 1}$ converges weakly to $J_{\gamma B}x$;
\item[(vi)] $(v^k)_{k\geq 1}$ converges weakly to $J_{\gamma B}x$;
\item[(vii)] if $A$ or $B$ is uniformly monotone, then $(y^k)_{k\geq 1}$ and $(v^k)_{k\geq 1}$ converge strongly to the unique 
point in $\zer(A+B)$.                    
\end{itemize} 
\end{theorem}

\begin{remark}\label{alpha1-neq0} According to \cite{b-c-h-inertial}, the condition $\alpha_1=0$ can be replaced with 
the assumption $w^0=w^1$ without altering the conclusion  of the above theorem.  
\end{remark}

\begin{remark}\label{expr<1} Let us mention that in the hypotheses of the above theorem we have 
$$0<\frac{\delta-\alpha\Big[\alpha(1+\alpha)+\alpha\delta+\sigma\Big]}{\delta\Big[1+\alpha(1+\alpha)+\alpha\delta+\sigma\Big]}<1.$$

\noindent Conversely, for a fixed $\ol\alpha\in(0,1)$, if we chose $\alpha> 0$ and $\sigma> 0$ such that 

$$\ol\alpha\big(1+\alpha(1+\alpha)+\sigma\big)+\alpha^2+2\alpha\sqrt{\ol\alpha}\sqrt{\alpha(1+\alpha)+\sigma}<1,$$

then $$\frac{\delta-\alpha\Big[\alpha(1+\alpha)+\alpha\delta+\sigma\Big]}{\delta\Big[1+\alpha(1+\alpha)+\alpha\delta+\sigma\Big]}=\ol\alpha,$$
for all $\delta\in\{\delta_1,\delta_2\}$, where 

$$\delta_{1,2}=\frac{1-\alpha^2-\ol\alpha\big(1+\alpha(1+\alpha)+\sigma\big)\pm
\sqrt{\left(1-\alpha^2-\ol\alpha\big(1+\alpha(1+\alpha)+\sigma\big)\right)^2-4\ol\alpha\alpha^2(\alpha(1+\alpha)+\sigma)}}{2\alpha\ol\alpha}.$$
 
\end{remark}

\section{Inertial ADMM}\label{IADMM}

In this section we present the main result of the paper, which consists in the formulation of an inertial ADMM algorithm for a primal-dual pair of convex optimization problems and in the investigation of its convergence properties. We start by describing the setting in which we work.

\begin{problem}\label{admm-p1} Let ${\cal H}$ and ${\cal G}$ be real Hilbert spaces, $f\in\Gamma({\cal H})$, $g\in\Gamma(\cal G)$ and $L:{\cal H}\rightarrow{\cal G}$ a linear continuous operator. We aim to solve the convex optimization problem 
\begin{equation}\label{prim} (P) \ \ \ \inf_{x\in{\cal H}}\{f(x)+g(Lx)\}
\end{equation}
together with its Fenchel-type dual problem 
\begin{equation}\label{dual} (D) \ \ \ \sup_{v\in{\cal G}}\{-f^*(-L^*v)-g^*(v)\}.
\end{equation}
\end{problem}

Denoting by $v(P)$ and $v(D)$ the optimal objective values of the two problems, respectively, the situation $v(P)\geq v(D)$, called in the literature {\it weak duality}, always holds. In case a regularity condition is fulfilled one can guarantee equality for the optimal objective values and existence of optimal solutions to the dual. For the readers convenience, we discuss some regularity conditions which are suitable in this context. One of the weakest regularity  conditions of interiority-type is the {\it Attouch-Br\'{e}zis condition}, which reads
\begin{equation}\label{reg-cond} 0\in\sqri(\dom g-L(\dom f)).
\end{equation}
Here, for $S\subseteq {\cal G}$ a convex set,  we denote by
$$\sqri S:=\{x\in S:\cup_{\lambda>0}\lambda(S-x) \ \mbox{is a closed linear subspace of} \ {\cal G}\}$$
its \textit{strong quasi-relative interior}. Notice that we always have $\inte S\subseteq\sqri S$ 
(in general this inclusion may be strict). If ${\cal G}$ is finite-dimensional, then $\sqri S$ coincides with 
$\ri S$, the relative interior of $S$, which is the interior of $S$ with respect to its affine hull. In this case, condition 
\eqref{reg-cond} holds if there exists $x'\in\ri(\dom f)$ such that $Lx'\in\ri(\dom g)$. Considering again the 
infinite dimensional setting, we remark that condition 
\eqref{reg-cond} is fulfilled, if for example $g$ is continuous at $x'\in\dom f\cap L^{-1}(\dom g)$. 
Let us mention that, if \eqref{reg-cond} holds, then we have {\it strong duality}, which means that $v(P)=v(D)$ and $(D)$ has an  optimal solution. 

Moreover, the optimality conditions for the primal-dual pair of optimization problems \eqref{prim}-\eqref{dual} read
\begin{equation}\label{opt-cond} -L^*v\in\partial f(x) \mbox{ and } v\in\partial g(Lx). 
\end{equation}

More precisely, if $(P)$ has an optimal solution $x\in{\cal H}$ and the regularity condition \eqref{reg-cond} is fulfilled, 
then there exists $v\in{\cal G}$, an optimal solution to $(D)$, such that  \eqref{opt-cond} holds. Conversely, if the pair 
$(x,v)\in{\cal H}\times{\cal G}$ satisfies relation \eqref{opt-cond}, then $x$ is an optimal solution to $(P)$ and  $v$ is an optimal solution to $(D)$. For further considerations concerning duality we invite the reader to consult \cite{bo-van, b-hab, bauschke-book, bot-csetnek, EkTem, simons, Zal-carte}.

Let us mention some conditions ensuring that $(P)$ has an optimal solution. 
Suppose that $(P)$ is feasible, which means that its optimal objective
value is not identical $+\infty$. The existence of optimal solutions to $(P)$ is guaranteed if,
for instance, $f$ is coercive (that is $\lim_{\|x\|\rightarrow\infty}f(x)=+\infty$)
and $g$ is bounded from below. Indeed, under these circumstances, the objective function of
$(P)$ is coercive and the statement follows via \cite[Corollary 11.15]{bauschke-book}. 
On the other hand, when $f$ is strongly convex, then the objective function of
$(P)$ is strongly convex, too, thus $(P)$ has a unique optimal solution (see \cite[Corollary 11.16]{bauschke-book}).

Let us introduce now the inertial ADMM algorithm. 

\begin{algorithm}\label{inertial-admm-alg} Chose $y^0,y^1,z^0,z^1\in{\cal G}$, $\gamma > 0$, $(\alpha_k)_{k\geq 1}$ nondecreasing with 
$0\leq\alpha_k\leq\alpha<1$ for every $k\geq 2$, $(\lambda_k)_{k\geq 1}$ and $\lambda, \sigma, \delta >0$ such that 
\begin{equation*}
\delta>\frac{\alpha^2(1+\alpha)+\alpha\sigma}{1-\alpha^2} \ \mbox{and} \ 0<\lambda\leq\lambda_k\leq 2\cdot\frac{\delta-\alpha\Big[\alpha(1+\alpha)+\alpha\delta+\sigma\Big]}{\delta\Big[1+\alpha(1+\alpha)+\alpha\delta+\sigma\Big]} \
\forall k \geq 2.  
\end{equation*}
Suppose that either $\alpha_2=0$ or $\lambda_1=\alpha_1=0$. Further, for all $k\geq 1$ set 
\begin{eqnarray} x^{k+1} & = &
\argmin_{x\in{\cal H}} \left\{f(x)+\left\langle y^k-\alpha_k(y^k-y^{k-1})-\gamma\alpha_k(z^k-z^{k-1}),Lx \right\rangle + \frac{\gamma}{2}\|Lx-z^k\|^2\right\} \label{inertial-x}\\
 \ol z^{k+1} & = & \alpha_{k+1}\lambda_k (Lx^{k+1}-z^k)+
\frac{(1-\lambda_k)\alpha_k\alpha_{k+1}}{\gamma}\left(y^k-y^{k-1}+\gamma(z^k-z^{k-1})\right)\label{inertial-zbar}\\
z^{k+1} & = &\argmin_{z\in{\cal G}}  \left\{g(z+\ol z^{k+1})+\left\langle -y^k-(1-\lambda_k)\alpha_k\left(y^k-y^{k-1}+\gamma(z^k-z^{k-1})\right),z\right\rangle \right .\nonumber \\
& & \qquad \qquad \qquad \quad   \quad  \ \  \left . +\frac{\gamma}{2}\|z-\lambda_kLx^{k+1}-(1-\lambda_k)z^k\|^2\right\} \label{inertial-z}\\
y^{k+1} & = & y^k+\gamma \left(\lambda_kLx^{k+1}+(1-\lambda_k)z^k-z^{k+1}\right)+
(1-\lambda_k)\alpha_k\left(y^k-y^{k-1}+\gamma(z^k-z^{k-1})\right).\label{inertial-y}
\end{eqnarray}
\end{algorithm}

\begin{remark} In order to ensure that the sequence $(x^k)_{k\geq 2}$ is uniquely determined we assume that the operator $L$ 
satisfies the hypothesis 

\begin{equation}\label{h} (H) \ \ \ \exists \theta>0 \mbox{ such that } \|Lx\|\geq\theta\|x\| \mbox{ for all }x\in{\cal H}. 
\end{equation}
This condition guarantees that the objective function in \eqref{inertial-x} is strongly convex, hence 
$(x^k)_{k\geq 2}$ is well defined (see \cite[Corollary 11.16]{bauschke-book}). Let us mention that $(H)$ will be used also in  
the proof of the convergence statements of the algorithm. Notice that if $L$ injective and $\ran L^*$ is closed, then 
$(H)$ holds (see \cite[Fact 2.19]{bauschke-book}). Moreover, $(H)$ implies that $L$ is injective. We conclude that 
in case $\ran L^*$ is closed, $(H)$ is equivalent to $L$ injective. In finite dimensional spaces, namely, if ${\cal H}=\R^n$ and 
${\cal G}=\R^m$, with $m\geq n\geq 1$, hypothesis $(H)$ is nothing else than saying that $L$ has full column rank, 
which is a condition widely used in the literature for proving the convergence of the ADMM algorithm.
\end{remark}

\begin{remark}\label{prox} Notice that the objective function of \eqref{inertial-z} is strongly convex, hence the sequence $(z^k)_{k\in\N}$ is well defined as well. Moreover, it can be expressed with the help of the proximal point operator of $g$ for every $k \geq 1$ as follows:
$$z^{k+1}=-\ol z^{k+1}+\prox\nolimits_{\gamma^{-1}g}\left(\ol z^{k+1}+\lambda_kLx^{k+1}+(1-\lambda_k)z^k+\frac{1}{\gamma}y^k+
\frac{(1-\lambda_k)\alpha_k}{\gamma}\left(y^k-y^{k-1}+\gamma(z^k-z^{k-1})\right)\right).$$
This is is general not the case for \eqref{inertial-x}, due to the presence of the operator $L$ in the $x$-argument. Nevertheless, in case ${\cal H}={\cal G}$ 
and $L$ is the identity operator on ${\cal H}$, relation \eqref{inertial-x} can be expressed via the proximal point operator of $f$ for every $k \geq 1$ as follows: 
$$x^{k+1}=\prox\nolimits_{\gamma^{-1}f}\left(z^k-\frac{1}{\gamma}y^k+\frac{\alpha_k}{\gamma}(y^k-y^{k-1})+\alpha_k(z^k-z^{k-1})\right).$$
\end{remark}

\begin{remark} Let us consider the case $\alpha_k=0$ for all $k\geq 1$. Then the iterative scheme becomes for every $k \geq 1$
\begin{eqnarray} x^{k+1} & = &
\argmin_{x\in{\cal H}}\left\{f(x)+\langle y^k,Lx \rangle
+\frac{\gamma}{2}\|Lx-z^k\|^2\right\} \label{inertial-x-lambda}\\
 z^{k+1} & = & \argmin_{z\in{\cal G}}\left\{g(z)+\langle -y^k,z\rangle
+\frac{\gamma}{2}\|z-\lambda_kLx^{k+1}-(1-\lambda_k)z^k\|^2\right\} \label{inertial-z-lambda}\\
y^{k+1} & =  & y^k+\gamma\left(\lambda_kLx^{k+1}+(1-\lambda_k)z^k-z^{k+1}\right) \label{inertial-y-lambda},
\end{eqnarray}
which is the error-free case of the classical ADMM algorithm as presented and investigated in \cite{ecb}. Here $(\lambda_k)_{k\geq 1}$ can be regarded as a sequence of relaxation parameters. If one takes further $\lambda_k=1$ for all $k\geq 1$, one has the 
classical ADMM algorithm (see for example \cite{bpcpe})
\begin{eqnarray} x^{k+1} & = &
\argmin_{x\in{\cal H}}\left\{f(x)+\langle y^k,Lx \rangle
+\frac{\gamma}{2}\|Lx-z^k\|^2\right\} \label{inertial-x-lambda1}\\
 z^{k+1} & = & \argmin_{z\in{\cal G}}\left\{g(z)+\langle -y^k,z\rangle
+\frac{\gamma}{2}\|z-Lx^{k+1}\|^2\right\} \label{inertial-z-lambda1}\\
y^{k+1} & =  & y^k+\gamma\left(Lx^{k+1}-z^{k+1}\right) \label{inertial-y-lambda1}.
\end{eqnarray}
\end{remark}

We are now in position to state the main result of the paper. 

\begin{theorem}\label{inertial-admm-th} 
In Problem \ref{admm-p1} suppose that $(P)$ has an optimal solution, the regularity 
condition \eqref{reg-cond} is fulfilled, the hypothesis $(H)$ concerning the operator $L$ holds and consider the sequences 
generated by Algorithm \ref{inertial-admm-alg}. Then there exists $(\ol x,\ol v)\in{\cal H}\times{\cal G}$ satisfying the 
optimality conditions \eqref{opt-cond}, hence, $\ol x$ is an optimal solution to $(P)$, $\ol v$ is an optimal solution  to $(D)$ and $v(P)=v(D)$, such that the following statements are true: 
\begin{enumerate}
 \item[(i)] $(x^k)_{k\geq 2}$ converges weakly to $\ol x$;
 \item[(ii)] $(\ol z^k)_{k\geq 2}$ converges strongly to $0$;
 \item[(iii)] $(z^k)_{k\in\N}$ converges weakly to $L\ol x$;
 \item[(iv)] $(Lx^{k+1}-z^k)_{k\geq 2}$ converges strongly to $0$;
 \item[(v)] $(y^k)_{k\in\N}$ converges weakly to $\ol v$;
 \item[(vi)] if $g^*$ is uniformly convex, then $(y^k)_{k\in\N}$ converges strongly to the unique optimal solution of $(D)$;
 \item[(vii)] $\lim_{k\rightarrow+\infty}(f(x^{k+1})+g(z^k+\ol z^k))=v(P)=v(D)=
              \lim_{k\rightarrow+\infty}(-f^*(-L^*v^k)-g^*(y^k))$, where the sequence $(v^k)_{k\geq 1}$ is defined by 
              \begin{equation}\label{inertial-v} v^k=y^k-\gamma z^k +\gamma Lx^{k+1}-\alpha_k\left(y^k-y^{k-1}+
               \gamma(z^k-z^{k-1})\right) \ \forall k\geq 1,
              \end{equation}
              and $(v^k)_{k\geq 1}$ converges weakly to $\ol v$.
\end{enumerate}
\end{theorem}

\begin{remark} Let us mention that the function $g^*$ is uniformly convex, if $g^*$ is $\beta$-strongly convex for $\beta > 0$.
Moreover, according to \cite[Theorem 18.15]{bauschke-book}, $g^*$ is $\beta$-strongly convex if and only if 
$g$ is Fr\'{e}chet-differentiable  and $\nabla g$ is $\beta^{-1}$-Lipschitzian.
\end{remark}

\begin{proof} 
We introduce the sequence $(w^k)_{k\in\N}$ by 
\begin{equation}\label{inertial-w} w^k=y^k+\gamma z^k \ \forall k\in\N.
\end{equation}

We intend to prove that the sequences $(y^k)_{k \in \N}, (v^k)_{k \geq 1}, (w^k)_{k \in \N}$ are nothing else than the ones generated by inertial Douglas-Rachford algorithm presented in 
Theorem \ref{inertial-DR} for the maximal monotone operators 
\begin{equation}\label{A,B} A:=\partial (f^*\circ (-L^*)) \mbox{ and } B:=\partial g^*.
\end{equation}
Notice that the hypotheses of the theorem ensure that there exist a pair $(x,v)\in{\cal H}\times{\cal G}$ satisfying the optimality conditions 
\eqref{opt-cond}, from which one easily derives that $\zer(A+B)\neq\emptyset$.

We fix $k\geq 1$. We obtain from \eqref{inertial-z} that 
$$0\in \partial g(z^{k+1}+\ol z^{k+1})-y^k-(1-\lambda_k)\alpha_k\left(y^k-y^{k-1}+\gamma(z^k-z^{k-1})\right)
+\gamma\left(z^{k+1}-\lambda_kLx^{k+1}-(1-\lambda_k)z^k\right),$$
hence due to \eqref{inertial-y} 
\begin{equation}\label{opt-cond-g} y^{k+1}\in \partial g(z^{k+1}+\ol z^{k+1}).
\end{equation}

From here we deduce $z^{k+1}+\ol z^{k+1}\in\partial g^*(y^{k+1})=By^{k+1}$, hence 
\begin{equation}\label{j-gb1}y^{k+1}=J_{\gamma B}(\gamma z^{k+1}+\gamma \ol z^{k+1}+y^{k+1})=J_{\gamma B}(w^{k+1}+\gamma \ol z^{k+1}).\end{equation}

By \eqref{inertial-y} we have 
\begin{equation}\label{cond-y} y^{k+1}=y^k+\gamma z^k-\gamma z^{k+1}+\gamma \lambda_k(Lx^{k+1}-z^k)+(1-\lambda_k)\alpha_k
\left(y^k-y^{k-1}+\gamma(z^k-z^{k-1})\right),
\end{equation}
thus \begin{equation}\label{ol z} \gamma \ol z^{k+1}=\alpha_{k+1}(y^{k+1}-y^k-\gamma z^k+\gamma z^{k+1})=\alpha_{k+1}(w^{k+1}-w^k).
     \end{equation}
     
From \eqref{j-gb1} and \eqref{ol z} we obtain 

\begin{equation}\label{dr1} y^{k+1}=J_{\gamma B}[w^{k+1}+\alpha_{k+1}(w^{k+1}-w^k)].
\end{equation}

Further, from \eqref{inertial-x} we get
$$0\in\partial f(x^{k+1})+L^*\left(y^k-\alpha_k(y^k-y^{k-1})-\gamma\alpha_k(z^k-z^{k-1})\right)+\gamma L^*(Lx^{k+1}-z^k),$$
which by \eqref{inertial-v} gives \begin{equation}\label{opt-cond-f} -L^*v^k\in\partial f(x^{k+1}). 
                                  \end{equation}
We derive $x^{k+1}\in\partial f^*(-L^*v^k)$, hence 
$$-Lx^{k+1}\in-L\partial f^*(-L^*v^k)\subseteq \partial (f^*\circ (-L^*))(v^k)=Av^k,$$
which leads to 
\begin{equation}\label{j-ga1} v^k=J_{\gamma A}(v^k-\gamma Lx^{k+1}).
\end{equation}

Taking into account \eqref{inertial-w} and \eqref{inertial-v} we have 
\begin{align*} v^k-\gamma Lx^{k+1}
& = y^k-\gamma z^k-\alpha_k\left(y^k-y^{k-1}+\gamma(z^k-z^{k-1})\right)\nonumber\\
& = 2y^k-w^k-\alpha_k\left(y^k-y^{k-1}+\gamma(z^k-z^{k-1})\right)\nonumber\\
& = 2y^k-w^k-\alpha_k(w^k-w^{k-1})\nonumber\\
\end{align*}
and from \eqref{j-ga1} we get \begin{equation}\label{dr2} v^k=J_{\gamma A}[2y^k-w^k-\alpha_k(w^k-w^{k-1})].
                              \end{equation}
Finally, from \eqref{cond-y}, \eqref{inertial-w} and \eqref{inertial-v} we derive 
\begin{align*} w^{k+1}
& = w^k+\alpha_k\left(y^k-y^{k-1}+\gamma(z^k-z^{k-1})\right)+\lambda_k\left(\gamma(Lx^{k+1}-z^k)-\alpha_k
\left(y^k-y^{k-1}+\gamma(z^k-z^{k-1})\right)\right)\nonumber\\
& = w^k+\alpha_k(w^k-w^{k-1})+\lambda_k(v^k-y^k),\nonumber\\
\end{align*}
hence \begin{equation}\label{dr3} w^{k+1}=w^k+\alpha_k(w^k-w^{k-1})+\lambda_k(v^k-y^k).
      \end{equation}

In conclusion, for all $k\geq 2$ we have (see \eqref{dr1}, \eqref{dr2} and \eqref{dr3}) $$\left\{
\begin{array}{ll}
y^k=J_{\gamma B}[w^k+\alpha_k(w^k-w^{k-1})]\\
v^k=J_{\gamma A}[2y^k-w^k-\alpha_k(w^k-w^{k-1})]\\
w^{k+1}=w^k+\alpha_k(w^k-w^{k-1})+\lambda_k(v^k-y^k), 
\end{array}\right.$$
which is the inertial Douglas-Rachford scheme from Theorem \ref{inertial-DR}. 

Notice that the relation $\alpha_2=0$ from Algorithm \ref{inertial-admm-alg} corresponds to the situation 
when in Theorem \ref{inertial-DR} the vectors $w^1,w^2$ can be chosen arbitrarily in ${\cal G}$, while the condition $\lambda_1=\alpha_1=0$ ensures that $w^1=w^2$, which is the situation mentioned in Remark \ref{alpha1-neq0}. Indeed, 
in case $\alpha_2\neq 0$ and $\lambda_1=\alpha_1=0$, from \eqref{inertial-zbar} we get $\ol z^2=0$, hence, by \eqref{ol z}, $w^1=w^2$. 

According to Theorem \ref{inertial-DR}, there exists $\ol w\in{\cal G}$ such that
\begin{equation}\label{th-dr1} w^k\rightharpoonup\ol w \ \mbox{as} \ k \rightarrow +\infty
\end{equation}
\begin{equation}\label{th-dr2} w^{k+1}-w^k\rightarrow 0 \ \mbox{as} \ k \rightarrow +\infty
\end{equation}
\begin{equation}\label{th-dr3} y^k-v^k\rightarrow 0 \ \mbox{as} \ k \rightarrow +\infty
\end{equation}
\begin{equation}\label{th-dr4} y^k\rightharpoonup J_{\gamma B}\ol w \ \mbox{as} \ k \rightarrow +\infty
\end{equation}
\begin{equation}\label{th-dr5} v^k\rightharpoonup J_{\gamma B}\ol w \ \mbox{as} \ k \rightarrow +\infty.
\end{equation}

From \eqref{ol z} and \eqref{th-dr2} we derive that \begin{equation}\label{ii}\ol z^k\rightarrow 0 \ \mbox{as} \ k \rightarrow +\infty. \end{equation} Further, 
by \eqref{inertial-v}, \eqref{th-dr2} and \eqref{th-dr3} we obtain  \begin{equation}\label{lx-z}Lx^{k+1}-z^k\rightarrow 0 \ \mbox{as} \ k \rightarrow +\infty.\end{equation}
Moreover, from \eqref{inertial-w}, \eqref{th-dr1} and \eqref{th-dr4} we get 
\begin{equation}\label{z}z^k\rightharpoonup\frac{1}{\gamma}(\ol w-J_{\gamma B}\ol w) \ \mbox{as} \ k \rightarrow +\infty.\end{equation} 
We deduce from \eqref{lx-z} that  \begin{equation}\label{lx}Lx^k\rightharpoonup\frac{1}{\gamma}(\ol w-J_{\gamma B}\ol w) \ \mbox{as} \ k \rightarrow +\infty.\end{equation}
Now, using the hypothesis $(H)$, we easily derive that $(x^k)_{k\geq 2}$ is bounded, thus, due to \eqref{lx}, it possesses at most 
one weak cluster point. As a consequence, $(x^k)_{k\geq 2}$ is weakly convergent (see \cite[Lemma 2.38]{bauschke-book}), hence 
there exists $\ol x\in{\cal H}$ such that \begin{equation}\label{x}x^k\rightharpoonup\ol x \ \mbox{as} \ k \rightarrow +\infty.\end{equation}
From \eqref{lx-z}, \eqref{z} and \eqref{x} we also have 
\begin{equation}\label{z2}z^k\rightharpoonup L\ol x \ \mbox{as} \ k \rightarrow +\infty\end{equation} and
\begin{equation}\label{lx-w}L\ol x=\frac{1}{\gamma}(\ol w-J_{\gamma B}\ol w).\end{equation}
Using the notation $\ol v=J_{\gamma B}\ol w$, we prove that the pair $(\ol x,\ol v)\in{\cal H}\times{\cal G}$ satisfies 
the optimality conditions \eqref{opt-cond}. To this end, observe that, due to \eqref{opt-cond-f} and \eqref{opt-cond-g}, we have 
\begin{equation}\label{prod-sp} (-L^*v^{k+1}+L^*y^{k+1},z^{k+1}+\ol z^{k+1}-Lx^{k+2})\in(\partial f\times B+S)(x^{k+2},y^{k+1}) \ \forall k\geq 1,
\end{equation}
where $S:{\cal H}\times{\cal G}\rightarrow{\cal H}\times{\cal G}$ is defined by
$$S(x,y)=(L^*y,-Lx) \ \forall (x,y)\in{\cal H}\times{\cal G}.$$
Since $S$ is monotone and continuous, it is maximally monotone (see \cite[Corollary 20.25]{bauschke-book}). Further, 
$\partial f\times B$ is also maximally monotone (see \cite[Proposition 20.23]{bauschke-book}) and since $S$ has full 
domain, the sum $\partial f\times B+S$ is maximally monotone, too (see \cite[Corollary 24.4]{bauschke-book}). Since the graph 
of a maximally monotone operator is sequentially closed in the weak-strong topology (see \cite[Proposition 20.33(ii)]{bauschke-book}), 
by taking the limits in \eqref{prod-sp} and using \eqref{th-dr3}, \eqref{lx-z}, \eqref{ii}, \eqref{x} and \eqref{th-dr4} 
we obtain $$(0,0)\in(\partial f\times B+S)(\ol x,\ol v).$$ One can easily show that the latter means the pair $(\ol x,\ol v)$ 
satisfies the optimality conditions \eqref{opt-cond}.

The statements (i)-(v) follow now from \eqref{x}, \eqref{ii}, \eqref{z2}, \eqref{lx-z} and \eqref{th-dr4}. Further, (vi) follows from  Theorem \ref{inertial-DR}(vii). 

We are going to prove now statement (vii). Notice that $f$ and $g$ are weak lower semicontinuous (since $f$ and $g$ are convex) and therefore, by (i), (ii) and (iii) 
we get \begin{align}\label{liminf} \liminf_{k\rightarrow+\infty} (f(x^{k+1})+g(z^k+\ol z^k)) 
                     & \geq \liminf_{k\rightarrow+\infty} f(x^{k+1}) + \liminf_{k\rightarrow+\infty} g(z^k+\ol z^k)\nonumber\\
                     & \geq f(\ol x) + g(L\ol x)=v(P).
       \end{align}
Further, from \eqref{opt-cond-f} we derive the inequality
\begin{equation}\label{f-x} f(\ol x)\geq f(x^{k+1})+\langle-L^*v^k,\ol x-x^{k+1}\rangle \ \forall k\geq 1  \end{equation}
and from \eqref{opt-cond-g} 
\begin{equation}\label{g-Lx} g(L\ol x)\geq g(z^k+\ol z^k)+\langle y^k,L\ol x-z^k-\ol z^k\rangle \ \forall k\geq 2. \end{equation}

Summing up the last two inequalities we get $$v(P)\geq f(x^{k+1})+g(z^k+\ol z^k)+\langle-v^k,L\ol x-Lx^{k+1}\rangle
+\langle y^k,L\ol x-z^k-\ol z^k\rangle \ \forall k\geq 2,$$
hence
$$f(x^{k+1})+g(z^k+\ol z^k)\leq v(P)+\langle v^k-y^k,L\ol x-Lx^{k+1}\rangle+\langle y^k,-Lx^{k+1}+z^k+\ol z^k\rangle \ \forall k\geq 2.$$
Taking into account \eqref{th-dr3}, (ii), (iii), (iv) and (v) we obtain
\begin{equation}\label{limsup}\limsup_{k\rightarrow+\infty} (f(x^{k+1})+g(z^k+\ol z^k))\leq v(P).  \end{equation}
Combining \eqref{liminf} and \eqref{limsup} we get the first part of the statement. 

Again by \eqref{opt-cond-f} and \eqref{opt-cond-g} we have (see \cite[Proposition 16.9]{bauschke-book}) 
\begin{equation}\label{f-x-f*} f(x^{k+1})+f^*(-L^*v^k)=\langle x^{k+1},-L^*v^k\rangle \ \forall k\geq 1  \end{equation}
and  
\begin{equation}\label{g-Lx-g*} g(z^k+\ol z^k)+g^*(y^k)=\langle y^k,z^k+\ol z^k\rangle \ \forall k\geq 2 . \end{equation}

Adding these relations we derive for every $k\geq 2$
$$-f^*(-L^*v^k)-g^*(y^k)=f(x^{k+1})+g(z^k+\ol z^k)+\langle v^k-y^k,Lx^{k+1}\rangle+\langle y^k,Lx^{k+1}-z^k-\ol z^k\rangle.$$
Finally, by \eqref{th-dr3}, (ii), (iii), (iv), (v) and the first part of (vii) we obtain 
$$\lim_{k\rightarrow+\infty}(-f^*(-L^*v^k)-g^*(y^k))=v(P)=v(D)$$ and the proof is complete. 
\end{proof}

\begin{remark} When working in finite dimensional spaces, there is no need for the construction considered in 
\eqref{prod-sp}, since in this case one can simply take the limits in \eqref{opt-cond-f} and \eqref{opt-cond-g} in order to 
conclude that $(\ol x,\ol v)$ satisfies the optimality conditions \eqref{opt-cond}. In infinite dimensional spaces this 
{\it naive} procedure does not work anymore, since in \eqref{opt-cond-f} and \eqref{opt-cond-g} we have only weak convergence 
for the sequences involved (we refer to \cite[Example 20.34]{bauschke-book} for an example of a maximally monotone operator whose 
graph is not sequentially closed in the weak-weak topology).
\end{remark}

\begin{remark}\label{f*-g*} Let us notice that the conclusion of Theorem \ref{inertial-admm-th}(vi) remains true if the uniform convexity 
of $g^*$ is replaced by the assumptions that $f^*$ is $\beta$-strongly convex, with $\beta >0$ and 
\begin{equation}\label{h*} (H^*) \ \ \ \exists \theta_*>0 \mbox{ such that } \|L^*v\|\geq\theta_*\|v\| \mbox{ for all }v\in{\cal G}. 
\end{equation}
 
Indeed, under these conditions one can prove that the composition $f^*\circ(-L^*)$ is $\beta\theta_*^2$-strongly convex, hence 
the operator $A$ (see \eqref{A,B}) is strongly monotone and the conclusion follows from Theorem \ref{inertial-DR}(vii). 
\end{remark}

\section{The minimization of a finite sum of convex functions}\label{IADMM-sum}

The aim of this section is to derive from Theorem \ref{inertial-admm-th} via the product space approach iterative schemes and corresponding convergence statements when solving the optimization problem which assumes the minimization of the finite sum of proper, convex and lower semicontinuous functions and its Fenchel-type dual. The goal is to evaluate each of the functions arising in the objective separately in the algorithmic scheme. 

\begin{problem}\label{admm-p2} Let ${\cal H}$ be a real Hilbert space, $m\geq 2$ a positive integer and 
$f_i\in\Gamma({\cal H})$ for $i=1,...,m$.. We aim to solve the convex optimization problem 
\begin{equation}\label{prim-sum} (P^{\sum}) \ \ \ \inf_{x\in{\cal H}}\left\{\sum_{i=1}^mf_i(x)\right\}
\end{equation}
together with its Fenchel-type dual problem 
\begin{equation}\label{dual-sum} (D^{\sum}) \ \ \ \sup_{\substack{v_i\in{\cal H},i=1,...,m\\\sum_{i=1}^m v_i=0}}\left\{-\sum_{i=1}^m f_i^*(v_i)\right\}.
\end{equation}
\end{problem}

One of the regularity conditions which guarantees strong duality in this situation is (see \cite{b-hab}): 
\begin{equation}\label{reg-cond-sum} 0\in\sqri\Big(\Pi_{i=1}^m\dom f_i-\{(x,...,x):x\in{\cal H}\}\Big).
\end{equation}

According to \cite[Remark 2.5]{b-hab}, this condition is fulfilled if there exists $x'\in \Pi_{i=1}^m\dom f_i$ such that 
$m-1$ functions $f_i$ are continuous at $x'$. In finite dimensional spaces, condition \eqref{reg-cond-sum} holds 
if $\cap_{i=1}^m\ri(\dom f_i)\neq\emptyset$. Also, let us mention that in case $m=2$, the regularity condition \eqref{reg-cond-sum} 
is equivalent to $0\in\sqri(\dom f_1-\dom f_2)$ (see \cite[Remark 2.5]{b-hab}). 

The optimality conditions for the primal-dual pair of optimization problems \eqref{prim-sum}-\eqref{dual-sum} read
\begin{equation}\label{opt-cond-sum} v_i\in\partial f_i(x) \ i=1,...,m \mbox{ and } \sum_{i=1}^mv_i=0. 
\end{equation}

More precisely, if $(P^{\sum})$ has an optimal solution $x\in{\cal H}$ and the regularity condition \eqref{reg-cond-sum} is fulfilled, 
then there exists $(v_1,...,v_m)\in{\cal H}^m$, an optimal solution to $(D^{\sum})$, such that \eqref{opt-cond-sum} holds. Conversely, if 
$(x,v_1,...,v_m)\in{\cal H}\times{\cal H}^m$ satisfies relation \eqref{opt-cond-sum}, then $x$ is an optimal solution to $(P^{\sum})$ and 
$(v_1,...,v_m)$ is an optimal solution to $(D^{\sum})$. 

Let us mention some conditions ensuring that $(P^{\sum})$ has an optimal solution. 
Suppose that $(P^{\sum})$ is feasible, which means that its optimal objective
value is not identical $+\infty$. The existence of optimal solutions of $(P^{\sum})$ is guaranteed if
for instance, one of the functions $f_i$  is coercive and the remaining ones are 
bounded from below. Indeed, under these circumstances, the objective function of
$(P^{\sum})$ is coercive and the statement follows via \cite[Corollary 11.15]{bauschke-book}. 
On the other hand, if one of the functions $f_i$ is strongly convex, then the objective function of
$(P^{\sum})$ is strongly convex, too, thus $(P^{\sum})$ has a unique optimal solution (see \cite[Corollary 11.16]{bauschke-book}).

We derive in the following two inertial ADMM algorithms for solving \eqref{prim-sum}-\eqref{dual-sum}. To this end we reformulate 
Problem \ref{admm-p2} as Problem \ref{admm-p1} in the product space ${\cal H}^m$ endowed with the inner product and associated norm defined by 
$$\langle x,u\rangle_{{\cal H}^m}=\sum_{i=1}^{m}\langle x_i,u_i\rangle_{\cal {H}} \mbox{ and }\|x\|_{{\cal H}^m}=\left(\sum_{i=1}^{m} \|x_i\|^2_{\cal H}\right)^{1/2}$$ 
for $x=(x_i)_{1\leq i\leq m}, u=(u_i)_{1\leq i\leq m} \in {\cal H}^m$, respectively,
where $\langle \cdot,\cdot \rangle_{\cal H}$ and $\|\cdot\|_{\cal H}$ denote the inner product and norm on ${\cal H}$, respectively. 

By using the notation $$C=\{(x,...,x):x\in{\cal H}\},$$
one can easily rewrite \eqref{prim-sum} as
\begin{equation}\label{prim-sum-prim}  \ \inf_{(x_1,...,x_m)\in{\cal H}^m}\Big\{f(x_1,...,x_m)+\delta_C(x_1,...,x_m)\Big\},
\end{equation}
where $f:{\cal H}^m\rightarrow\B$ is defined by
$$f(x_1,...,x_m)=\sum_{i=1}^mf_i(x_i) \ \forall (x_1,...,x_m)\in{\cal H}^m.$$
This corresponds to the optimization problem \eqref{prim} with $g=\delta_C:{\cal H}^m\rightarrow\B$ and $L$ is the identity operator 
on ${\cal H}^m$. Notice that the Fenchel dual problem \eqref{dual} of \eqref{prim-sum-prim} becomes $(D^{\sum})$, the regularity condition 
\eqref{reg-cond} is equivalent to \eqref{reg-cond-sum} and the optimality conditions \eqref{opt-cond} are nothing else than the ones in
\eqref{opt-cond-sum}. Moreover, $(x_1,...,x_m)\in{\cal H}^m$ is an optimal solution of \eqref{prim-sum-prim} if and only if 
$x_i=x$, $i=1,..,m,$ and $x\in{\cal H}$ is an optimal solution to $(P^{\sum})$, while, $(v_1,...,v_m)\in{\cal H}^m$ is an optimal 
solution to the dual of \eqref{prim-sum-prim} if and only if $(v_1,...,v_m)\in{\cal H}^m$ is an optimal 
solution to $(D^{\sum})$. This shows that we are in the context of Problem \ref{admm-p1}. 

Writing Algorithm \ref{inertial-admm-alg} in this setting we get for every $k\geq 1$ the following iterative scheme 
\begin{eqnarray} x^{k+1} & = &
\argmin_{x\in{\cal H}^m}\left\{f(x)+\left\langle y^k-\alpha_k(y^k-y^{k-1})-\gamma\alpha_k(z^k-z^{k-1}),x \right\rangle_{{\cal H}^m}
+\frac{\gamma}{2}\|x-z^k\|_{{\cal H}^m}^2\right\} \label{inertial-x-sum1'} \\
 \ol z^{k+1} & = & \alpha_{k+1}\lambda_k(x^{k+1}-z^k)+
\frac{(1-\lambda_k)\alpha_k\alpha_{k+1}}{\gamma}\left(y^k-y^{k-1}+\gamma(z^k-z^{k-1})\right) \label{inertial-zbar-sum1'}\\
 z^{k+1} & = &
\argmin_{z\in{{\cal H}^m}}\left\{g(z+\ol z^{k+1})+\left\langle -y^k-(1-\lambda_k)\alpha_k\left(y^k-y^{k-1}+\gamma(z^k-z^{k-1})\right),z\right\rangle_{{\cal H}^m} \right . \nonumber\\
& & \qquad \qquad \qquad \qquad \ \ \ \left . +\frac{\gamma}{2}\left\|z-\lambda_kx^{k+1}-(1-\lambda_k)z^k\right\|_{{\cal H}^m}^2\right\} \label{inertial-z-sum1'}\\
 y^{k+1} & = & y^k+\gamma\left(\lambda_kx^{k+1}+(1-\lambda_k)z^k-z^{k+1}\right)+
(1-\lambda_k)\alpha_k\left(y^k-y^{k-1}+\gamma(z^k-z^{k-1})\right), \label{inertial-y-sum1'}
\end{eqnarray}
where $x^{k+1}=(x_i^{k+1})_{1\leq i\leq m}$, $\ol z^{k+1}=(\ol z_i^{k+1})_{1\leq i\leq m}$, $z^{k+1}=(z_i^{k+1})_{1\leq i\leq m}$, 
$y^{k+1}=(y_i^{k+1})_{1\leq i\leq m}$, $x=(x_i)_{1\leq i\leq m}$ and $z=(z_i)_{1\leq i\leq m}$. 

We give in the following an explicit form of this algorithm. Due to the definition of $f$, relation \eqref{inertial-x-sum1'} is nothing else 
than: \begin{equation}\label{inertial-x-sum1-i} x_i^{k+1}=
\argmin_{x\in{\cal H}}\left\{f_i(x)+\left\langle y_i^k-\alpha_k(y_i^k-y_i^{k-1})-\gamma\alpha_k(z_i^k-z_i^{k-1}),x \right\rangle
+\frac{\gamma}{2}\|x-z_i^k\|^2\right\}, \ i=1,...,m.
\end{equation}

Further, from \eqref{inertial-z-sum1'} and \eqref{opt-cond-g} we derive $$z^{k+1}+\ol z^{k+1}\in C$$ and $$y^{k+1}\in C^{\perp},$$
hence there exists $(u^k)_{k\geq 2}\in{\cal H}$ such that for every $k\geq 1$ it holds \begin{equation}\label{z-C} z_i^{k+1}+\ol z_i^{k+1}=u^{k+1}, \  i=1,...,m \end{equation}
and \begin{equation}\label{y-C-perp}\sum_{i=1}^my_i^{k+1}=0.\end{equation}

If we suppose that $\sum_{i=1}^my_i^k=0$ for every $k\geq 0$, then from \eqref{inertial-y-sum1'} we derive 
\begin{equation}\label{sum1-z} \sum_{i=1}^m z_i^{k+1}= \lambda_k\sum_{i=1}^mx_i^{k+1}+(1-\lambda_k)\sum_{i=1}^m z_i^k +(1-\lambda_k)
\alpha_k\sum_{i=1}^m(z_i^k-z_i^{k-1}) \ \forall k\geq 1.
\end{equation}
From this, \eqref{z-C} and \eqref{inertial-zbar-sum1'}, we get
\begin{equation}\label{u} u^{k+1}=\frac{\lambda_k(1+\alpha_{k+1})}{m}\sum_{i=1}^mx_i^{k+1}+\frac{1-\alpha_{k+1}\lambda_k-\lambda_k}{m}\sum_{i=1}^mz_i^k
+\frac{\alpha_k(1-\lambda_k)(1+\alpha_{k+1})}{m}\sum_{i=1}^m(z_i^k-z_i^{k-1}) \ \forall k\geq 1.\end{equation}
Conversely, if for a fixed $k\geq 1$ we suppose that $\sum_{i=1}^my_i^{k-1}=\sum_{i=1}^my_i^k=0$, then from \eqref{inertial-y-sum1'}, \eqref{z-C} and 
\eqref{u} we have $\sum_{i=1}^my_i^{k+1}=0$. 

All together, we derive the following algorithm and corresponding convergence theorem (notice that for the statement in (vii) we use also Remark \ref{f*-g*}). 

\begin{algorithm}\label{inertial-admm-alg-sum1} Chose $y_i^0,y_i^1,z_i^0,z_i^1\in{\cal H}$, $i=1,...,m$, such that $\sum_{i=1}^my_i^0=\sum_{i=1}^my_i^1=0$, $\gamma > 0$, $(\alpha_k)_{k\geq 1}$ nondecreasing with 
$0\leq\alpha_k\leq\alpha<1$ for every $k\geq 2$, $(\lambda_k)_{k\geq 1}$ and $\lambda, \sigma, \delta >0$ such that 
\begin{equation*}
\delta>\frac{\alpha^2(1+\alpha)+\alpha\sigma}{1-\alpha^2} \ \mbox{and} \ 0<\lambda\leq\lambda_k\leq 2\cdot\frac{\delta-\alpha\Big[\alpha(1+\alpha)+\alpha\delta+\sigma\Big]}{\delta\Big[1+\alpha(1+\alpha)+\alpha\delta+\sigma\Big]} \
\forall k \geq 2.  
\end{equation*}
Suppose that either $\alpha_2=0$ or $\lambda_1=\alpha_1=0$. Further, for every $k\geq 1$ set 
\begin{eqnarray}
\label{inertial-x-sum1} x_i^{k+1} & = & \argmin_{x\in{\cal H}}\left\{f_i(x)+\left\langle y_i^k-\alpha_k(y_i^k-y_i^{k-1})-\gamma\alpha_k(z_i^k-z_i^{k-1}),x \right\rangle
+\frac{\gamma}{2}\|x-z_i^k\|^2\right\},\ i=1,...,m\\
\label{inertial-zbar-sum1} \ol z_i^{k+1} & = & \alpha_{k+1}\lambda_k(x_i^{k+1}-z_i^k)+
\frac{(1-\lambda_k)\alpha_k\alpha_{k+1}}{\gamma}\left(y_i^k-y_i^{k-1}+\gamma(z_i^k-z_i^{k-1})\right), \  i=1,...,m\\
\label{inertial-u-sum1} u^{k+1} & = & \frac{\lambda_k(1+\alpha_{k+1})}{m}\sum_{i=1}^mx_i^{k+1}+\frac{1-\alpha_{k+1}\lambda_k-\lambda_k}{m}\sum_{i=1}^mz_i^k
+\frac{\alpha_k(1-\lambda_k)(1+\alpha_{k+1})}{m}\sum_{i=1}^m(z_i^k-z_i^{k-1}) \\ 
\label{inertial-z-sum1} z_i^{k+1} & = & u^{k+1}-\ol z_i^{k+1}, \ i=1,...,m\\
\label{inertial-y-sum1} y_i^{k+1} & = & y_i^k+\gamma\left(\lambda_kx_i^{k+1}+(1-\lambda_k)z_i^k-z_i^{k+1}\right)+
(1-\lambda_k)\alpha_k\left(y_i^k-y_i^{k-1}+\gamma(z_i^k-z_i^{k-1})\right), \ i=1,...,m.
\end{eqnarray}
\end{algorithm}

\begin{theorem}\label{inertial-admm-th-sum1} In Problem \ref{admm-p2} suppose that $(P^{\sum})$ has an optimal solution, the regularity 
condition \eqref{reg-cond-sum} is fulfilled and consider the sequences 
generated by Algorithm \ref{inertial-admm-alg-sum1}. Then there exists $(\ol x,\ol v_1,...,\ol v_m)\in{\cal H}\times{\cal H}^m$ satisfying the 
optimality conditions \eqref{opt-cond-sum}, hence $\ol x$ is an optimal solution to $(P^{\sum})$, $(\ol v_1,...,\ol v_m)$ is an optimal solution 
to $(D^{\sum})$ and $v(P^{\sum})=v(D^{\sum})$, such that the following statements are true: 
\begin{enumerate}
 \item[(i)] $(x_i^k)_{k\geq 2}$ converges weakly to $\ol x$, $i=1,...,m$;
 \item[(ii)] $(\ol z_i^k)_{k\geq 2}$ converges strongly to $0$, $i=1,...,m$;
 \item[(iii)] $(z_i^k)_{k\in\N}$ converges weakly to $\ol x$, $i=1,...,m$;
 \item[(iv)] $(x_i^{k+1}-z_i^k)_{k\geq 2}$ converges strongly to $0$, $i=1,...,m$;
 \item[(v)] $(u^k)_{k\geq 2}$ converges weakly to $\ol x$; 
 \item[(vi)] $(y_i^k)_{k\in\N}$ converges weakly to $\ol v_i$, $i=1,...,m$;
 \item[(vii)] if $f_i^*$ is strongly convex for every $i=1,...,m$, then $((y_1^k)_{k\in\N},...,(y_m^k)_{k\in\N})$ converges strongly to the unique optimal
             solution to $(D^{\sum})$;
 \item[(viii)] $\lim_{k\rightarrow+\infty}\left(\sum_{i=1}^mf_i(x_i^{k+1})\right)=v(P^{\sum})=v(D^{\sum})=
              \lim_{k\rightarrow+\infty}\left(-\sum_{i=1}^mf_i^*(-v_i^k)\right)$, where for every $i=1,...,m$, the sequence $(v_i^k)_{k\geq 1}$ is defined by 
              \begin{equation}\label{inertial-v-sum1} v_i^k=y_i^k-\gamma z_i^k +\gamma x_i^{k+1}-\alpha_k\left(y_i^k-y_i^{k-1}+
               \gamma(z_i^k-z_i^{k-1})\right) \ \forall k\geq 1,
              \end{equation}
              and $(v_i^k)_{k\geq 1}$ converges weakly to $\ol v_i$.
\end{enumerate}
\end{theorem}

\begin{remark} If we take $\lambda_k=1$ for every $k\geq 1$, then $\ol z_i^{k+1}=\alpha_{k+1}(x_i^{k+1}-z_i^k), i=1,...,m,$ and (see also relation \eqref{sum1-z}) 
$u^{k+1}=\frac{1+\alpha_{k+1}}{m}\sum_{i=1}^mx_i^{k+1}-\frac{\alpha_{k+1}}{m}\sum_{i=1}^mx_i^k$, hence 
the iterative scheme \eqref{inertial-x-sum1} - \eqref{inertial-y-sum1} can be simplified to  
\begin{eqnarray}
\label{inertial-x-sum1-lambda1} x_i^{k+1} & = & \argmin_{x\in{\cal H}}\left\{f_i(x)+\left\langle y_i^k-\alpha_k(y_i^k-y_i^{k-1})-\gamma\alpha_k(z_i^k-z_i^{k-1}),x \right\rangle
+\frac{\gamma}{2}\|x-z_i^k\|^2\right\}, \ i=1,...,m\\
\label{inertial-z-sum1-lambda1} z_i^{k+1} & = & \frac{1+\alpha_{k+1}}{m}\sum_{j=1}^mx_j^{k+1}-\frac{\alpha_{k+1}}{m}\sum_{j=1}^mx_j^k- \alpha_{k+1}(x_i^{k+1}-z_i^k),\ i=1,...,m\\
\label{inertial-y-sum1-lambda1} y_i^{k+1} & = & y_i^k+\gamma\left(x_i^{k+1}-z_i^{k+1}\right), \ i=1,...,m.
\end{eqnarray}
If, moreover, $\alpha_k=0$ for every $k\geq 1$, \eqref{inertial-x-sum1-lambda1} - \eqref{inertial-y-sum1-lambda1} becomes 
\begin{eqnarray}
\label{inertial-x-sum1-lambda1-alpha0} x_i^{k+1} & = & \argmin_{x\in{\cal H}}\left\{f_i(x)+\langle y_i^k,x\rangle
+\frac{\gamma}{2}\left\|x-\frac{1}{m}\sum_{j=1}^mx_j^k\right\|^2\right\}, \  i=1,...,m\\
\label{inertial-y-sum1-lambda1-alpha0} y_i^{k+1} & = & y_i^k+\gamma\left(x_i^{k+1}-\frac{1}{m}\sum_{j=1}^mx_j^{k+1}\right), \ i=1,...,m,
\end{eqnarray}
which is the ADMM algorithm as considered in \cite[page 50]{bpcpe}. 
\end{remark}

By interchanging the roles of $f$ and $g$ in \eqref{prim-sum-prim} we obtain another inertial ADMM-type algorithm with corresponding convergence statement. 

\begin{algorithm}\label{inertial-admm-alg-sum2} Chose $y_i^0,y_i^1,z_i^0,z_i^1\in{\cal H}$, $i=1,...,m$, $\gamma > 0$, $(\alpha_k)_{k\geq 1}$ nondecreasing with 
$0\leq\alpha_k\leq\alpha<1$ for every $k\geq 2$, $(\lambda_k)_{k\geq 1}$ and $\lambda, \sigma, \delta >0$ such that 
\begin{equation*}
\delta>\frac{\alpha^2(1+\alpha)+\alpha\sigma}{1-\alpha^2} \ \mbox{and} \ 0<\lambda\leq\lambda_k\leq 2\cdot\frac{\delta-\alpha\Big[\alpha(1+\alpha)+\alpha\delta+\sigma\Big]}{\delta\Big[1+\alpha(1+\alpha)+\alpha\delta+\sigma\Big]} \
\forall k \geq 2.  
\end{equation*}
Suppose that either $\alpha_2=0$ or $\lambda_1=\alpha_1=0$. Further, for every $k\geq 1$ set 
\begin{eqnarray}
\label{inertial-x-sum2} x^{k+1} & = & \frac{1}{m}\sum_{i=1}^mz_i^k-\frac{1}{m\gamma}\sum_{i=1}^my_i^k+
\frac{\alpha_k}{m\gamma}\sum_{i=1}^m\left(y_i^k-y_i^{k-1}+\gamma(z_i^k-z_i^{k-1})\right)\\
\label{inertial-zbar-sum2} \ol z_i^{k+1} & = & \alpha_{k+1}\lambda_k(x^{k+1}-z_i^k)+
\frac{(1-\lambda_k)\alpha_k\alpha_{k+1}}{\gamma}\left(y_i^k-y_i^{k-1}+\gamma(z_i^k-z_i^{k-1})\right), \  i=1,...,m\\
\label{inertial-z-sum2} z_i^{k+1} & = & \argmin_{z\in{\cal H}}\left\{f_i(z+\ol z_i^{k+1})+\left\langle -y_i^k-(1-\lambda_k)\alpha_k\left(y_i^k-y_i^{k-1}+\gamma(z_i^k-z_i^{k-1})\right),z\right\rangle \right.\\
& & \qquad  \qquad  \qquad  \qquad \ \ \ \left . +\frac{\gamma}{2}\|z-\lambda_kx^{k+1}-(1-\lambda_k)z_i^k\|^2\right\}, \ i=1,...,m\\
\label{inertial-y-sum2} y_i^{k+1} & = & y_i^k+\gamma\left(\lambda_kx^{k+1}+(1-\lambda_k)z_i^k-z_i^{k+1}\right)+(1-\lambda_k)\alpha_k\left(y_i^k-y_i^{k-1}+\gamma(z_i^k-z_i^{k-1})\right), \  i=1,...,m.
\end{eqnarray}
\end{algorithm}

\begin{theorem}\label{inertial-admm-th-sum2} In Problem \ref{admm-p2} suppose that $(P^{\sum})$ has an optimal solution, the regularity 
condition \eqref{reg-cond-sum} is fulfilled and consider the sequences 
generated by Algorithm \ref{inertial-admm-alg-sum2}. Then there exists $(\ol x,\ol v_1,...,\ol v_m)\in{\cal H}\times{\cal H}^m$ satisfying the 
optimality conditions \eqref{opt-cond-sum}, hence $\ol x$ is an optimal solution to $(P^{\sum})$, $(\ol v_1,...,\ol v_m)$ is an optimal solution 
to $(D^{\sum})$ and $v(P^{\sum})=v(D^{\sum})$, such that the following statements are true: 
\begin{enumerate}
 \item[(i)] $(x^k)_{k\geq 2}$ converges weakly to $\ol x$;
 \item[(ii)] $(\ol z_i^k)_{k\geq 2}$ converges strongly to $0$, $i=1,...,m$;
 \item[(iii)] $(z_i^k)_{k\in\N}$ converges weakly to $\ol x$, $i=1,...,m$;
 \item[(iv)] $(x_i^{k+1}-z_i^k)_{k\geq 2}$ converges strongly to $0$,  $i=1,...,m$;
  \item[(v)] $(y_i^k)_{k\in\N}$ converges weakly to $\ol v_i$, $i=1,...,m$;
 \item[(vi)] if $f_i^*$ is strongly convex for $i=1,...,m$, then $((y_1^k)_{k\in\N},...,(y_m^k)_{k\in\N})$ converges strongly to the unique optimal
             solution to $(D^{\sum})$;
 \item[(vii)] $\lim_{k\rightarrow+\infty}\left(\sum_{i=1}^mf_i(z_i^k+\ol z_i^k)\right)=v(P^{\sum})=v(D^{\sum})=
              \lim_{k\rightarrow+\infty}\left(-\sum_{i=1}^mf_i^*(-y_i^k)\right)$.
\end{enumerate}
\end{theorem}

\begin{remark} Notice that relation \eqref{inertial-x-sum2} is derived from $v^k\in C^{\perp}$ (see \eqref{opt-cond-f}), where  
for any $i=1,...,m$ the sequence $(v_i^k)_{k\geq 1}$ is defined by 
              \begin{equation}\label{inertial-v-sum2} v_i^k=y_i^k-\gamma z_i^k +\gamma x^{k+1}-\alpha_k\left(y_i^k-y_i^{k-1}+
               \gamma(z_i^k-z_i^{k-1})\right) \ \forall k\geq 1
              \end{equation}
              and $(v_i^k)_{k\geq 1}$ converges weakly to $\ol v_i$.
\end{remark}

\end{document}